\documentclass[11pt]{article}

\usepackage{amssymb,latexsym}
\usepackage{amsmath}
\usepackage{citesort}
\usepackage{eucal}

\title{Deformation Quantization of Hermitian Vector Bundles}

\author{{\bf     
         Henrique Bursztyn\thanks{henrique@math.berkeley.edu}
         \thanks{Research supported by a fellowship from CNPq, Grant
           200481/96-7.} 
        } \\[0.5cm]
        Department of Mathematics\\
        UC Berkeley\\
        94720 Berkeley, CA, USA
        \\[1cm]
        {\bf
         Stefan
         Waldmann\thanks{Stefan.Waldmann@ulb.ac.be}
         \thanks{Research 
         supported by the Communaut\'e fran\c caise de Belgique,
         through an Action de Recherche Concert\'ee de la Direction de la
         Recherche Scientifique.}
        } \\[0.5cm]
        D{\'{e}}partement de Math{\'{e}}matique \\
        Universit{\'e} Libre de Bruxelles \\
        Campus Plaine, C. P. 218 \\
        Boulevard du Triomphe \\
        B-1050 Bruxelles \\
        Belgique 
       }

\date{September 2000}

\newcommand{\im} {{\mathrm i}}

\newcommand{\cc} [1]     {\overline {{#1}}}

\newcommand{\id}         {{\mathsf {id}}}
\newcommand{\tr}         {\mathop{{\mathsf {tr}}}}

\newcommand{\End}        {\mathop{{\mathsf {End}}}}
\newcommand{\ring}[1]    {{\mathsf {{#1}}}}

\newcommand{\SP} [1]     {{\left\langle {{#1}} \right\rangle}}
\newcommand{\Bounded}    {{\mathfrak B}}
\newcommand{\Unit}       {\mathsf {1}}

\newcommand{\Mloc}       {\mathcal{M}_{\mathrm {loc}}}
\newcommand{\Def}        {\mathsf{Def}_{\mathrm {loc}}}
\newcommand{\Right}      {{R}}

\newcommand{\qh} {\boldsymbol{h}}
\newcommand{\qR} {\boldsymbol{R}}
\newcommand{\qA} {\boldsymbol{\mathcal A}}
\newcommand{\qB} {\boldsymbol{\mathcal B}}
\newcommand{\qE} {\boldsymbol{\mathcal E}}
\newcommand{\qP} {\boldsymbol{P}}
\newcommand{\qQ} {\boldsymbol{Q}}
\newcommand{\qU} {\boldsymbol{U}}
\newcommand{\qV} {\boldsymbol{V}}

\newcommand{\qx} {\boldsymbol{x}}
\newcommand{\qy} {\boldsymbol{y}}
\newcommand{\qT} {\boldsymbol{T}}
\newtheorem{lemma} {Lemma} [section]
\newtheorem{proposition} [lemma] {Proposition}
\newtheorem{theorem} [lemma] {Theorem}
\newtheorem{corollary} [lemma] {Corollary}
\newtheorem{definition}[lemma] {Definition}

\newenvironment{proof}{\small{\sc Proof:}}{{\hspace*{\fill} $\square$\\}}

\numberwithin{equation}{section}

\begin{document}

\maketitle

\begin{abstract}
Motivated by deformation quantization, we consider in this paper
$^*$-algebras $\mathcal A$ over rings $\ring C = \ring{R}(\im)$,
where $\ring R$ is an ordered ring and $\im^2 = -1$, and study
the deformation theory of projective
modules over these algebras carrying the additional structure of
a (positive) $\mathcal A$-valued inner product. For $A=C^\infty(M)$,
$M$ a manifold, these modules can be identified with Hermitian
vector bundles $E$ over $M$. We show that for a fixed
Hermitian star-product on $M$,  these modules can 
always be deformed in a unique way, up to (isometric) equivalence.
We observe that there is a natural bijection between the sets of
equivalence classes of local Hermitian deformations of $C^\infty(M)$
and $\Gamma^\infty(\End(E))$ and that the corresponding deformed
algebras are formally Morita equivalent, an algebraic generalization
of strong Morita equivalence of $C^*$-algebras. We also discuss the
semi-classical geometry arising from these deformations.
\end{abstract}
\newpage

\section{Introduction} \label{IntroSec}

The concept of formal deformation quantization was first introduced in
\cite{BFFLS78} and its goal is to construct quantum observable 
algebras by means of formal deformations (in the sense of Gerstenhaber
\cite{GS88}) of classical Poisson algebras, usually given by the
algebra of complex-valued smooth functions on a Poisson manifold. The
deformed associative algebra structures arising in this way are called
star-products. Their existence and classification have been
established through the joint effort of many authors (see 
\cite{DL83b,Fed94a,OMY91,Kon97b,NT95a,NT95b,BCG97,WX98} 
and \cite{Wei94,Ste98,Gutt2000} for surveys on the subject).
The aim of the present paper is to study a certain class of modules
over these star-product algebras.

In order to better understand how deformation quantization relates
to the usual formalism of quantum mechanics through operators on
Hilbert spaces, one is naturally led to consider representations of
star-product algebras. By using the natural order structure in the
ring $\mathbb{R}[[\lambda]]$ and considering Hermitian star-products
(i.e., star-products for which the pointwise complex conjugation of
functions is a $^*$-involution), a theory of $^*$-representations of
star-product algebras on formal pre-Hilbert spaces can be developed,
similar to the usual representation theory of $C^*$-algebras.
In fact, these ideas can be carried out in the more general setting of
$^*$-algebras over $\ring C = \ring{R}(\im)$, where $\ring R$ is an
ordered ring and $\ring C$ is its ring extension by $\im$, with 
$\im^2 = -1$. This category encompasses both star-product algebras and
operator algebras. This approach was started in \cite{BW98a} with a
formal version of the GNS construction and has been further
investigated in \cite{BW97b,BNW99a,BNW98a,Wal2000a}. In this context,
one can formulate the notion of formal Morita equivalence
\cite{BuWa99a,BuWa2000}, which is a generalization of Rieffel's notion
of strong Morita equivalence of $C^*$-algebras \cite{Rief82}. An
important role in this theory is played by finitely generated
projective modules over  unital algebras $\mathcal A$ carrying the
additional structure of an $\mathcal A$-valued inner product. For
$\mathcal A = C^\infty(M)$, $M$ an arbitrary manifold, these modules
can be naturally identified with (smooth sections of) smooth Hermitian
vector bundles over $M$. The main subject of the present paper is the
deformation quantization of projective inner product $\mathcal
A$-modules. We remark that other authors have introduced
quantization of vector bundles in different contexts, such as
geometric
quantization, quantum groups and non-commutative geometry
\cite{Haw98,Haw99,BrMj93,Connes94,Rief88b}.

This paper is organized as follows. In Section~\ref{StarAlgSect} 
we collect some basic facts about $^*$-algebras over ordered
rings. Section~\ref{DefProjSec} reviews some aspects of the  
deformation theory of $^*$-algebras and their projections.
In Section~\ref{InnerProdSec}, we discuss deformations of projective
modules with inner-products. The main result in this section is
Proposition~\ref{HerEProp}. We emphasize here the additional Hermitian
structure, since the result without inner product seems to be
well-known, although not stated in the literature in this form. In
Section~\ref{DefHerVBSec}, we specialize the discussion of
Section~\ref{InnerProdSec} to the case of Hermitian vector bundles 
$E \to M$, where $M$ is a Poisson manifold with a star-product. We
show existence and uniqueness of deformation quantization of these
objects, see Definition~\ref{DeqQuanDef} and
Theorem~\ref{MainThm}. Using this result, we also show that there is a
natural bijection between the sets of equivalence classes of local
Hermitian deformations of the $^*$-algebras $C^\infty(M)$ and
$\Gamma^\infty(\End(E))$ in Proposition~\ref{BijecProp}. This
generalizes the result in \cite{LecRog88}, where the the case of
trivial bundles is treated through a different approach. We show in
Section~\ref{FullProjSec} that corresponding deformations are
formally Morita equivalent
(Proposition~\ref{StrongFullDeformProp}). The important ingredient for
this result is the notion of a strongly full projection (see
Definition~\ref{StrongFullDef}). A discussion about the underlying
semi-classical geometry corresponding to the deformations just
mentioned (in the sense of \cite{RVW96}) is presented in
Section~\ref{SemiClassSec}.

We remark that techniques used here have been previously used 
by other authors, see e.g. \cite{EmmWe96,Fed96,Ros96}. The idea of
deforming projections is present in Fedosov's index theorems
\cite{Fed96} (see also \cite{NT95a,NT95b}) and in the star-product 
formulation of multicomponent WKB expansions by Weinstein and Emmrich
\cite{EmmWe96} (see also \cite{ER98}). It would be interesting to
explore the connections between the present paper and these topics.
We also believe this paper is  an important step in extending the
approach to WKB approximations taken in \cite{BNW99a,BW97b} to the
multicomponent case.

Finally, let us mention that various authors in particle physics 
have considered ``noncommutative gauge fields'', see
e.g. \cite{MSSW2000,JSW2000} and the references therein. Here the
fields can be viewed as sections of deformed vector bundles, but
mainly trivial bundles have been treated so far. The Hermitian
metrics (and their deformations) should play an important role in this
context since they typically provide the potential term in the
Lagrangians of physical models. This will be the subject of future
investigations.

\section{Algebraic preliminaries} 
\label{PrelimSec}


\subsection{$^*$-algebras over ordered rings}
\label{StarAlgSect}

We will recall here some basic definitions concerning $^*$-algebras
over ordered rings. Further details can be found in
\cite{BuWa99a,BuWa2000,BW98a,BuWa99b}. All algebras in this paper will
be assumed to be associative.

Let $\ring{R}$ be an ordered ring, i.e. an associative, commutative,
and unital ring with a subset $\ring{P} \subset \ring{R}$ so that
$\ring R$ is the disjoint union 
$\ring{R} = -\ring{P} \cup \{0\} \cup \ring{P}$ and 
$\ring{P} \cdot \ring{P} \subseteq \ring{P}$, 
$\ring{P} + \ring{P} \subseteq \ring{P}$.
An element $a \in \ring{P}$ is called \emph{positive} and we denote it
by $a > 0$. We remark that if $\ring R$ is ordered, then the ring
$\ring{R}[[\lambda]]$ of formal power series has a natural ordering
given by $\sum_{r=0}^\infty a_r \lambda^r >0$ if $a_{r_0}>0$,
where $r_0$ is the first index with nonvanishing coefficient. The
important examples of ordered rings in deformation quantization are
$\mathbb{R}$ and $\mathbb{R}[[\lambda]]$. We also define 
$\ring C = \ring{R}(\im)$ to be the quadratic ring extension of 
$\ring R$ by $\im$, with $\im ^2 = -1$. Complex conjugation 
$z \mapsto \overline z$ is defined in the usual way. In the following
we shall assume for simplicity that $\mathbb Q \subseteq \ring R$.

Let $\mathcal A$ be an associative algebra over $\ring C$, equipped
with a $^*$-involution, i.e. an involutive $\ring C$-antilinear
antiautomorphism $^*: \mathcal A \longrightarrow \mathcal A$. We call
$\mathcal A$ a \emph{$^*$-algebra} over $\ring C$. We define
\emph{Hermitian, unitary} and \emph{normal} elements in $\mathcal A$
in the usual way. A linear functional 
$\omega: \mathcal A \longrightarrow \ring C$ is called \emph{positive}
if $\omega(A^*A)\geq 0$ for all $A \in \mathcal A$. An element 
$A \in \mathcal A$ is then called \emph{positive} if 
$\omega(A) \geq 0$ for all positive linear functionals $\omega$. The
set of positive elements in $\mathcal A$ is denoted by 
$\mathcal A^+$. Note that all elements of the form 
$a_1 A_1^*A_1 + \cdots + a_n A_n^*A_n$, 
with $A_j \in \mathcal A$ and $a_j \geq 0$ in $\ring C$, are
positive. These elements are called \emph{algebraically positive} and they
form a set denoted by $\mathcal A^{++}$. 
As an example of such $^*$-algebras, consider 
$\mathcal A = C^\infty(M)$, the algebra of complex-valued smooth
functions on a manifold $M$ with $^*$-involution given by complex
conjugation. In this case, with these definitions, positive
functionals correspond to positive Borel measures with compact support
and positive elements in $\mathcal A$ are positive functions as
expected, see e.g.~\cite[App.~B]{BuWa99a}. We remark that all these
notions also make sense for star-product algebras
$(C^\infty(M)[[\lambda]],\star)$ by means of the order structure of
$\mathbb{R}[[\lambda]]$, see \cite{BW97b,BuWa99b}.

\subsection{Formal deformations of $^*$-algebras and projections}
\label{DefProjSec}

We shall discuss here some aspects of the formal deformation theory of
associative algebras and $^*$-algebras over $\ring C$. A good part of
it is well-known and the reader is refered to \cite{GS88,BuWa99b} for
more details.

Let $k$ be a commutative and unital ring and let $\mathcal A$ be a
$k$-algebra. A \emph{formal deformation} of $\mathcal A$ (in the sense
of Gerstenhaber, see e.g.~\cite{GS88}) is an associative
$k[[\lambda]]$-bilinear multiplication $\star$ on
$\mathcal{A}[[\lambda]]$ of the form  
\begin{equation}
    A \star A' = \sum_{r=0}^\infty C_r(A,A') \lambda^r, 
    \qquad A,A' \in \mathcal A,
\end{equation}
where each $C_r$ is a Hochshild $2$-cochain and 
$C_0:\mathcal{A} \times \mathcal{A} \longrightarrow \mathcal A$ is the
original product on $\mathcal A$. 
If $\mathcal A$ is unital, we require in addition that the unit
element $\Unit \in \mathcal A$ is still the unit element with respect
to $\star$. We denote the deformed algebra by
$\qA = (\mathcal A[[\lambda]], \star)$. We recall that two
deformations of $\mathcal A$, 
$\qA_1 = (\mathcal{A}[[\lambda]],\star_1)$ and 
$\qA_2 = (\mathcal{A}[[\lambda]],\star_2)$, are called
\emph{equivalent} if there exist $k$-linear maps 
$T_r: \mathcal{A}\longrightarrow \mathcal{A}$, $r \geq 1$, so that 
$\qT = 
\id + \sum_{r=1}^\infty T_r \lambda^r : \qA_1 \longrightarrow \qA_2$ 
satisfies
\begin{equation}
    A \star_1 A' = \qT^{-1}(\qT(A)\star_2 \qT(A')), 
    \qquad \forall A, A' \in \mathcal A[[\lambda]].
\end{equation}

Let us now consider $\mathcal A$ to be a $^*$-algebra 
over $\ring C = \ring{R}(\im)$, $\ring{R}$ ordered. (For simplicity,
whenever we refer to a $^*$-algebra $\mathcal A$, it will be
implicitly assumed that the underlying ring $k$ is $\ring C$.) A
formal deformation $\qA = (\mathcal{A}[[\lambda]],\star)$ of 
$\mathcal A$ is called \emph{Hermitian} if the natural extension of
the $^*$-involution in $\mathcal A$ to $\mathcal{A}[[\lambda]]$ is
still an involution with respect to  $\star$, i.e. 
$(A \star A')^*={A'}^*\star A^*$ for all $A, A' \in \mathcal A$.
In this paper, deformations of $^*$-algebras will be always assumed to
be Hermitian. We recall that if $\qA_1$, $\qA_2$ are two Hermitian
deformations of $\mathcal A$ which are equivalent then there actually
exists an equivalence $\qT$ satisfying, in addition,
$\qT(A^*)=\qT(A)^*$, for all $A \in \mathcal A$ (see
\cite[Prop.~5.6]{Neu99}). So equivalence transformations between Hermitian
deformations will be assumed to preserve the involution.

We observe that if $\mathcal A$ is a $k$-algebra and $\qA$ is a formal
deformation of $\mathcal A$, then $M_n(\qA)$ can be identified with
$M_n(\mathcal A)[[\lambda]]$ as a $k[[\lambda]]$-module and naturally
defines a deformation of $M_n(\mathcal A)$. Also note that if
$\mathcal A$ is a $^*$-algebra, we can define a $^*$-involution on
$M_n(\mathcal A)$ in the usual way and if $\qA$ is a Hermitian
deformation of $\mathcal A$, $M_n(\qA)$ naturally defines a Hermitian
deformation of $M_n(\mathcal A)$. 
\begin{lemma}\label{TrickLem}
Let $\qA$ be a Hermitian deformation of a unital $^*$-algebra
$\mathcal A$ over $\ring C$.
Let $L_0 \in M_n(\mathcal A)$ be invertible and let 
$\boldsymbol{S} = \sum_{r=0}^\infty S_r \lambda^r \in M_n(\qA)$ be
Hermitian with $S_0 = L_0^*L_0$.
Then there exist $L_r \in M_n(\mathcal A), r \geq 1$, such that 
$\boldsymbol{L} = \sum_{r=0}^\infty L_r \lambda^r$ satisfies 
$\boldsymbol{S} = \boldsymbol{L}^* \star \boldsymbol{L}$.
\end{lemma}
\begin{proof}
We define $\boldsymbol{L}$ recursively. Suppose 
$L_0,L_1, \ldots, L_{k-1} \in M_n(\mathcal A)$ are such that 
$\boldsymbol{L}_{k-1} 
= L_0+ L_1 \lambda + \ldots + L_{k-1}\lambda^{k-1}$ 
satisfies 
$\boldsymbol{S} -\boldsymbol{L}_{k-1}^* \star \boldsymbol{L}_{k-1} 
= b_k \lambda^k + o(\lambda^{k+1})$. Note that since $\boldsymbol{S}$
is Hermitian, so is $b_k$. We need to find $L_k$ so that 
$\boldsymbol{L}_k = \sum_{j=0}^k L_j \lambda^j$ 
satisfies 
$\boldsymbol{S} = \boldsymbol{L}_{k}^*\star \boldsymbol{L}_{k}$ up to
order $\lambda^{k+1}$. But this happens if and only if 
$L_k^*L_0 + L_0^*L_k = b_k$. Then $L_k = \frac{1}{2}(b_kL_0^{-1})^*$
is a solution. 
\end{proof}
\begin{corollary}
Let $\mathcal A$ be a unital $^*$-algebra over $\ring C$ and $\qA$ a
Hermitian deformation of $\mathcal A$. Then any unitary 
$U_0 \in M_n(\mathcal A)$ can be deformed into a unitary 
$\qU = \sum_{r=0}^\infty U_r \lambda^r \in M_n(\qA)$.
\end{corollary}

Let $\mathcal A$ be a $k$-algebra and let $P_0 \in M_n(\mathcal A)$ be
an idempotent, i.e. $P_0^2=P_0$. It is well-kown that if $\qA$ is any
formal deformation of $\mathcal A$, we can always deform $P_0$, that
is, find an idempotent $\qP \in M_n(\qA)$ so that 
$\qP = P_0 + o(\lambda)$ (see e.g.~\cite{Fed96,EmmWe96,GS88}). 
In particular, the explicit formula (e.g.~\cite[Eq.~(6.1.4)]{Fed96})
\begin{equation}
    \label{DeformedProjection}
    \qP = \frac{1}{2} + \left(P_0 - \frac{1}{2}\right) \star
    \frac{1}{\sqrt[\star]{1 + 4 (P_0 \star P_0 - P_0)}}
\end{equation}
shows that, in the case $\mathcal A$ is a $^*$-algebra and $\qA$ is a 
Hermitian deformation, $\qP$ can be chosen to be a
projection (i.e. a \emph{Hermitian} idempotent) if $P_0$ is a projection.
\begin{lemma}\label{ILem}
Let $\mathcal A$ be a $k$-algebra and suppose $P_0, Q_0 \in M_n(\mathcal A)$ 
are idempotents. Let $\qA$ be a deformation of $\mathcal A$ and 
$\qP=\sum_{r=0}^\infty P_r \lambda^r
, \qQ= \sum_{r=0}^\infty Q_r\lambda^r
\in M_n(\qA)$ be deformations of $P_0, Q_0$, respectively. Then the map
$I: P_0M_n(\mathcal A)Q_0 [[\lambda]] \longrightarrow 
\qP \star M_n(\qA)\star \qQ$ given by
\begin{equation}
    \label{IEq}
    I(P_0 L Q_0) 
    = \qP \star (P_0 L Q_0) \star \qQ,
    \qquad L \in M_n(\mathcal A)[[\lambda]],
\end{equation}
is a $k[[\lambda]]$-module isomorphism.
\end{lemma}
\begin{proof}
The $k[[\lambda]]$-linearity and the
injectivity of $I$ are obvious since $\qP$ and $\qQ$ are deformations
of $P_0$ and $Q_0$. To prove surjectivity, let 
$\boldsymbol{L} = \sum_{r=0}^\infty L_r\lambda^r \in 
\qP \star M_n (\qA) \star\qQ$
be given. Then $\qP\star \boldsymbol{L} \star \qQ = 
\boldsymbol{L}$ whence $P_0 L_0 Q_0 =
L_0$. Thus defining $S_0 := L_0 \in P_0 M_n(\mathcal A)Q_0$, 
we have $I(S_0) = \boldsymbol{L}$ up to order $\lambda^0$. 
Since $I(S_0) - \boldsymbol{L} \in \qP \star M_n (\qA) \star\qQ$
starts with order $\lambda$, we can repeat the argument to find a 
$S_1 \in P_0 M_n(\mathcal A)Q_0$ such that $I(S_0 + \lambda S_1)$
coincides with $\boldsymbol{L}$ up to order $\lambda$. Then a simple induction
proves that $I$ is onto.
\end{proof}

Keeping the notation as in Lemma \ref{ILem}, we denote the deformed
product of $M_n(\qA)$ by 
$L \star S = \sum_{r=0}^\infty C_r(L,S) \lambda ^r$,
$L, S \in M_n(\mathcal A)$. Note that if we consider 
$I:P_0M_n(\mathcal A)Q_0[[\lambda]]
\longrightarrow M_n(\qA)=M_n(\mathcal A)[[\lambda]]$, we can write
$I=\sum_{r=0}^\infty I_r \lambda^r$ and a simple computation shows that
\begin{equation}
    \label{IrEq}
    I_r(B)= \sum_{i+j+k+m = r}C_m(C_k(P_i,B),Q_j), 
    \qquad \textrm{for} \quad B \in P_0 M_n(\mathcal A)Q_0.
\end{equation}
We note that $I$ is just a deformation of the natural inclusion
$P_0M_n(\mathcal A)Q_0 \hookrightarrow M_n(\mathcal A)$.

Let us consider a $k$-algebra $\mathcal A$, $P_0 \in M_n(\mathcal A)$
an idempotent, $\qA=(\mathcal{A}[[\lambda]], \star)$ a formal
deformation of $\mathcal A$ and $\qP$ a deformation of $P_0$. It is
clear that $P_0M_n(\mathcal A)P_0$ and $\qP \star M_n(\qA) \star \qP$
are algebras, which are unital if $\mathcal A$ is unital. Note that if
$\mathcal A$ is a $^*$-algebra, $P_0 \in M_n(\mathcal A)$ is a
projection, and $\qA$ is a Hermitian deformation of $\mathcal A$, with
$\qP$ a projection deforming $P_0$, then $P_0M_n(\mathcal A)P_0$ and
$\qP \star M_n(\qA) \star \qP$ are in fact $^*$-algebras. From
Lemma~\ref{ILem} we find 
\begin{corollary}\label{PAPCor}
The map 
$I:P_0M_n(\mathcal A)P_0[[\lambda]] \longrightarrow \qP\star M_n(\qA)
\star \qP$ 
defined in (\ref{IEq}) induces a formal deformation of
$P_0M_n(\mathcal A)P_0$. Moreover, if $\mathcal A$ is a $^*$-algebra,
$P_0$ is a projection and $\qA$ is Hermitian, then $I$ induces a
Hermitian deformation of $P_0M_n(\mathcal A)P_0$.
\end{corollary}

\subsection{Deformations of projective inner-product
            $\mathcal A$-modules} 
\label{InnerProdSec}


Let $\mathcal A$ be a unital $k$-algebra and let $\mathcal E$ be a
module over $\mathcal A$.  We
will essentially restrict ourselves to \emph{right} modules, but the
reader will have no problem to adapt all the definitons and results to
come to left modules.

Let 
$\Right_0: \mathcal{E} \times \mathcal{A} \longrightarrow \mathcal E$
denote the right $\mathcal A$-action on $\mathcal E$, 
$\Right_0(x, A) = x \cdot A$ for 
$x \in \mathcal E$, $A \in \mathcal A$. 
Let $\qA = (\mathcal{A}[[\lambda]], \star)$ be a formal deformation of
$\mathcal A$ and suppose there exist $k$-bilinear maps 
$\Right_r:\mathcal E \times \mathcal A \longrightarrow \mathcal E$,
for $r \geq 1$, such that the map
\begin{equation}
    \label{RrEq}
    \qR = \sum_{r=0}^\infty \Right_r \lambda^r : 
    \mathcal{E}[[\lambda]] \times \qA 
    \longrightarrow \mathcal{E}[[\lambda]]
\end{equation}
makes $\mathcal{E}[[\lambda]]$ into a module over $\qA$. We will
denote $\qR(x,A) = x \bullet A$, for $x \in \mathcal E$, 
$A \in \mathcal A$. 
\begin{definition}\label{qEDef}
We call $\qE=(\mathcal{E}[[\lambda]], \bullet)$ a \emph{deformation}
of the (right) $\mathcal A$-module $\mathcal E$ corresponding to 
$\qA=(\mathcal{A}[[\lambda]], \star)$. Two deformations 
$\qE=(\mathcal{E}[[\lambda]], \bullet)$, 
$\qE'=(\mathcal{E}[[\lambda]], \bullet')$ are \emph{equivalent} if
there exists an $\qA$-module isomorphism
$\qT: \qE \longrightarrow \qE'$ of the form 
$\qT = \id + \sum_{r=1}^\infty T_r \lambda^r$, 
with $k$-linear maps $T_r:\mathcal E \longrightarrow \mathcal E$.
\end{definition}
Here we will be only interested in finitely generated projective modules
(f.g.p.m.).
\begin{proposition}\label{DefModProp}
Let $\mathcal A$ be a unital $k$-algebra and
$\qA=(\mathcal{A}[[\lambda]], \star)$ be a deformation of 
$\mathcal A$. Let $\mathcal E$ be a (right) f.g.p.m. over 
$\mathcal A$. Then there exists a deformation $\qE$
of $\mathcal E$ corresponding to $\qA$ so that $\qE$ is a f.g.p.m.
over $\qA$. Moreover, this deformation is unique up to equivalence and
hence, in particular, every deformation of $\mathcal E$ is finitely
generated and projective. 
\end{proposition}
\begin{proof}
For the existence, note that since $\mathcal E$ is f.g.p.m., it
follows that we can identify $\mathcal E = P_0 \mathcal{A}^n$, for
some $n \geq 1$ and $P_0 \in M_n(\mathcal A)$ idempotent. Let 
$\qP \in M_n(\qA)$ be an idempotent deforming $P_0$ and consider the
(right) f.g.p.m. over $\qA$ given by $\qP \star \qA^n$. By
Lemma~\ref{ILem} (choosing $Q_0$ to be $1$ in the upper right corner
and zero elsewhere), we can use the isomorphism 
$I: \mathcal{E}[[\lambda]] \longrightarrow \qP \star \qA^n$ to pull
this $\qA$-module structure back to $\mathcal{E}[[\lambda]]$, i.e.
$x \bullet A := I^{-1}(\qP\star x \star A)= x \cdot A + o(\lambda)$ 
for $x \in \mathcal E$ and $A \in \mathcal A$. 
So $\qE = (\mathcal{E}[[\lambda]],\bullet)$ is a f.g.p. deformation of 
$\mathcal E$. 

Now assume 
$\qE'=(\mathcal{E}[[\lambda]],\bullet')$ is another deformation of
$\mathcal E$. Let 
$\mathfrak{C}: \qE \longrightarrow \qE/{\lambda \qE}=\mathcal{E}$ and
$\mathfrak{C}': \qE' \longrightarrow \qE'/{\lambda \qE'}=\mathcal{E}$
be the natural projections, which are surjective $\qA$-module
homomorphisms. Then it follows by projectivity of $\qE$ that there
exists an $\qA$-module homomorphism $\qT: \qE \longrightarrow \qE'$
satisfying $\mathfrak{C}'\circ \qT = \mathfrak{C}$. Since 
$\qE = \qE' = \mathcal{E}[[\lambda]]$ as
$k[[\lambda]]$-modules, we can write 
$\qT=\sum_{r=0}^\infty T_r\lambda^r$ and it is readly seen that
$T_0=\id$. So $\qT$ is an equivalence. 
\end{proof}

Suppose now $\mathcal A$ is a $^*$-algebra. A (right) f.g.p.m. 
$\mathcal E$ over $\mathcal A$ is called an 
\emph{inner-product module} if it is equipped 
with a positive definite, $\mathcal A$-valued, and $\mathcal A$-right
linear inner product, i.e. with a map
$h_0:\mathcal{E} \times \mathcal{E} \longrightarrow \mathcal A$
satisfying 
$h_0(x, \alpha y + \beta z) = \alpha h_0(x,y) + \beta h_0(x,z)$,
$h_0(x,y)^*=h_0(y,x)$,
$h_0(x, y \cdot A) = h_0 (x, y) A$, 
$h_0(x,x) \in \mathcal A^+$, and $h_0(x,x)=0 \implies x=0$, where
$x,y,z \in \mathcal E$, $\alpha, \beta \in \ring C$, and 
$A \in \mathcal A$. 
A left inner-product module is defined analogously but $h_0$ is
required to be $\ring C$-linear and $\mathcal A$-left linear in the
first argument. The $\mathcal A$-right linear endomorphisms of
$\mathcal E$ are denoted by $\End_{\mathcal A}(\mathcal E)$ and the
inner product determines the subalgebra 
$\Bounded_{\mathcal A}(\mathcal E) 
\subseteq \End_{\mathcal A} (\mathcal E)$
of those endomorphisms which have an adjoint with respect to
$h_0$. Then $\Bounded_{\mathcal A} (\mathcal E)$ becomes a
$^*$-algebra over $\ring C$.

Let $\qA$ be a Hermitian deformation of $\mathcal A$.
Let $\qE = (\mathcal{E}[[\lambda]], \bullet)$ be a corresponding
deformation of $\mathcal E$ and suppose there exist
$h_r:\mathcal{E}\times \mathcal{E} \longrightarrow \mathcal A$ such
that 
\begin{equation} 
    \label{hrEq}
    \qh=\sum_{r=0}^\infty h_r \lambda^r
\end{equation}
defines a positive definite, $\qA$-valued, $\qA$-right linear inner
product on $\mathcal{E}[[\lambda]]$.
\begin{definition}\label{qHEDef}
We call $\qE=(\mathcal{E}[[\lambda]], \bullet, \qh)$ a \emph{Hermitian
deformation} of the inner-product module $(\mathcal E, h_0)$ corresponding to
$\qA$. Two Hermitian deformations $\qE=(\mathcal{E}[[\lambda]], \bullet, \qh)$,
$\qE=(\mathcal{E}[[\lambda]], \bullet', \qh')$ are called \emph{equivalent}
if there is an equivalence 
$\qT=\id + \sum_{r=0}^\infty T_r \lambda^r : \qE \longrightarrow \qE'$
(as in Definition \ref{qEDef}) satisfying 
$\qh'(\qT(x),\qT(y))=\qh(x,y)$, $x,y \in \mathcal E$. 
\end{definition}

Let $\mathcal A$ be a unital $^*$-algebra and let 
$P_0 \in M_n(\mathcal A)$ be a projection. Consider the f.g.p.m. over
$\mathcal A$ given by $\mathcal E = P_0 \mathcal{A}^n$. We observe
that $\mathcal E$ has a canonical $\mathcal A$-valued inner product
$h_0$, namely the restriction to $P_0 \mathcal{A}^n$ of the canonical
$\mathcal A$-valued inner product on the free-module $\mathcal{A}^n$
given by $\SP{x,y}=\sum_{i=1}^n x_i^* y_i$. 
Note that in this case 
$M_n(\mathcal A) = \End_{\mathcal A}(\mathcal{A}^n) =
\Bounded_{\mathcal A} (\mathcal A^n)$ with the above $^*$-involution,
the same holding for 
$P_0M_n(\mathcal A)P_0 = \End_{\mathcal A}(P_0 \mathcal{A}^n) =
\Bounded_{\mathcal A} (P_0 \mathcal A^n)$.
\begin{proposition} \label{HerEProp}
Let $\mathcal A$ be a unital $^*$-algebra and 
$P_0 \in M_n(\mathcal A)$ be a projection. Let $\qA$ be a Hermitian
deformation of $\mathcal A$ and consider the $\mathcal A$-module
$\mathcal E = P_0\mathcal{A}^n$, equipped with its canonical 
$\mathcal A$-valued inner product $h_0$. Then there exists a 
Hermitian deformation of $\mathcal E$ corresponding to $\qA$, which is
unique up to equivalence.
\end{proposition}
\begin{proof}
As in Proposition \ref{DefModProp}, we choose $\qP \in M_n(\qA)$,
a projection deforming $P_0$, and consider the $\qA$-module
$\qP \star \qA^n$, which we know to define a deformation 
$\qE = (\mathcal E[[\lambda]], \bullet)$ of $\mathcal E$. 
Let $\qh$ be the $\qA$-valued inner product on $\qE$ obtained from
$\SP{\cdot,\cdot}$ restricted to $\qP \star \qA^n$. A simple 
computation shows that $\qh$ is a deformation of $h_0$ and hence 
$\qE=(\qP \star \qA^n,\qh)$ is a Hermitian deformation of 
$(\mathcal E, h_0)$.

Let $\qE'=(\mathcal{E}[[\lambda]],\bullet',\qh')$ be another Hermitian
deformation of $(\mathcal E, h_0)$. By Proposition \ref{DefModProp},
we may assume that $\qE' = \qE$ as an $\qA$-module, with
some $\qA$-valued inner product $\qh'$ deforming $h_0$. Recall that
any $\qA$-valued inner product 
$\boldsymbol{\langle}\cdot, \cdot \boldsymbol{\rangle}'$ 
on the free $\qA$-module $\qA^n$ can be written as
$\boldsymbol{\langle}\cdot, \cdot \boldsymbol{\rangle}' = 
\boldsymbol{\langle}\cdot, \boldsymbol{H} \cdot \boldsymbol{\rangle}$
for some Hermitian element $\boldsymbol{H}\in M_n(\qA)$, where 
$\boldsymbol{\langle}\qx, \qy \boldsymbol{\rangle}= 
\sum_{i=1}^n \qx_i^* \star \qy_i$, $\qx,\qy \in \qA^n$.
Since $\qP \star \qA^n \subseteq \qA^n$ is projective, one can check
that the same holds for this submodule. Hence, there is a Hermitian
element $\boldsymbol{H} \in \End_{\qA}(\qE)$ so that 
$\qh'(\cdot, \cdot) = \qh(\cdot, \boldsymbol{H} \cdot)$. 
But since $\qh$ and $\qh'$ are deformations of $h_0$, we can write 
$\boldsymbol{H} = \sum_{r=0}^\infty H_r \lambda^r$ 
with $H_0 = \id$. It then follows from Lemma \ref{TrickLem} that we
can find $\qU = \id + \sum_{r=0}^\infty U_r \lambda^r$ so that 
$\boldsymbol{H} = \qU^*\star\qU$. It is then clear that 
$\qU : \qE' \longrightarrow \qE$ is the desired equivalence.
\end{proof}

Let us finally discuss the deformation of isometries. An element 
$V_0 \in \End_{\mathcal A} (\mathcal E)$, with 
$\mathcal E = P_0 \mathcal A^n$, is called an \emph{isometry} of 
$\mathcal E$ if $V_0$ is invertible and $h_0(V_0x, V_0y) = h_0(x,y)$
for all $x, y \in \mathcal E$. Clearly such a $V_0$ gives rise to a
unitary $\tilde{V_0} \in M_n (\mathcal A)$ with 
$\tilde{V_0} P_0 = P_0 \tilde{V_0}$ and every such unitary yields an
isometry of $\mathcal E$ by restriction. 
\begin{proposition}
    \label{IsometryProp}
    Let $\qA$ be a Hermitian deformation of $\mathcal A$ and 
    $\qE = (\mathcal E[[\lambda]], \bullet, h)$ a Hermitian
    deformation of $\mathcal E = P_0 \mathcal A^n$. Then for every
    isometry $V_0 \in \End_{\mathcal A} (\mathcal E)$ there exists a
    deformation $\qV = V_0 + \sum_{r=1}^\infty \lambda^r V_r$ of $V_0$
    into an isometry $\qV$ of $\qE$.
\end{proposition}
\begin{proof}
    Since all deformations of $\mathcal E$ are equivalent we choose a
    deformation $\qP$ of the projection $P_0$. Moreover, we choose a
    unitary $\tilde{\qU} \in M_n (\qA)$ with 
    $\tilde{\qU} = \sum_{r=0}^\infty \lambda^r \tilde U_r$ and
    $\tilde U_0 = \tilde V_0$ which is possible due to
    Lemma~\ref{TrickLem}. Now 
    $\qP' = \tilde{\qU} \star \qP \star \tilde{\qU}^*$ is
    again a projection with classical limit $P_0$. Thus it defines a
    deformation $(\mathcal E[[\lambda]], \bullet', \qh')$ which is
    equivalent to the first one by an equivalence transformation 
    $\qT= \id + \sum_{r=1}^\infty \lambda^r T_r$. Moreover,
    $\tilde{\qU}$ descends to a unitary $\qA$-right linear map 
    $\qU: (\mathcal E[[\lambda]], \bullet, \qh) \to
    (\mathcal E[[\lambda]], \bullet', \qh')$ with lowest order 
    $U_0 = V_0$. Then $\qV = \qT^{-1} \circ \qU$ is the desired
    deformation of $V_0$. 
\end{proof}

\section{Deformation quantization of Hermitian vector bundles} 


\subsection{Deformation quantization}
\label{DefHerVBSec}

Let $\mathcal A = C^\infty(M)$ be the algebra of complex-valued smooth
functions on a manifold $M$. This algebra has a natural
$^*$-involution given by complex conjugation. Let $E \to M$ be a
complex vector bundle over $M$, with fiber dimension $k\geq
1$. Consider $\mathcal E = \Gamma^\infty(E)$, the space of smooth
sections of $E$, equipped with its natural right $\mathcal A$-module
structure. We recall that $\mathcal E$ is a f.g.p.m. over 
$\mathcal A$ (see \cite[Ch.~I, Thm.~6.5]{Karoubi78},
noticing that the proof there works for any manifold due to
\cite[Lem.~2.7]{Munkres63}). Finally, let 
$\mathcal B = \End_{\mathcal A}(\mathcal E) = \Gamma^\infty(\End(E))$ 
be the complex algebra of smooth sections of the endomorphism bundle
$\End(E) \to M$. Note that if $E$ is Hermitian, i.e. equipped with a
Hermitian fiber metric $h_0$, then there is a corresponding positive
definite $\mathcal A$-valued, $\mathcal A$-right linear inner product
on $\mathcal E$, also denoted by $h_0$. This defines a $^*$-involution
on $\mathcal B$.

Let $\qA = (C^\infty(M)[[\lambda]], \star)$ be a deformation of
$\mathcal A$. Recall that a star-product $\star$, given by 
$f \star g = \sum_{r=0}^\infty C_r(f,g) \lambda^r$,  is called
\emph{local/differential/of Vey type} if each $C_r$ is
local/differential/differential of order $r$ in each argument.
\begin{definition}
\label{DeqQuanDef}
Let $\qA = (C^\infty(M)[[\lambda]], \star)$ be a deformation of
$\mathcal A$. A \emph{deformation quantization} of $E$ is a
deformation of $\mathcal E = \Gamma^\infty(E)$ in the sense of
Definition \ref{qEDef}. If $E$ is equipped with a Hermitian fiber
metric $h_0$ and $\qA$ is a Hemitian deformation of $\mathcal A$, then
a \emph{Hermitian deformation quantization} of $(E, h_0)$ is a
deformation of $(\mathcal E, h_0)$ as in Definition \ref{qHEDef}. 
A deformation is called \emph{local/differential/of Vey type} if the 
corresponding $\Right_r$ and $h_r$ (as in (\ref{RrEq}), (\ref{hrEq}))
are local/differential/differential of order r in each argument.
\end{definition}

Recall that any two Hermitian metrics on a complex vector bundle are
equivalent (see \cite[Ch.~I,Thm.~8.8]{Karoubi78}), 
hence we can identify $\mathcal E$
with $P_0\mathcal{A}^n$, for some $n \geq 1$ and some projection 
$P_0 \in M_n(\mathcal A)$, equipped with its canonical 
$\mathcal A$-valued inner product $h_0$. 
\begin{theorem}\label{MainThm}
Let $E$ be a complex (Hermitian) vector bundle over a Poisson manifold
$M$ and let $\qA=(C^\infty(M)[[\lambda]], \star)$ be a (Hermitian)
deformation quantization of $\mathcal A$. Then there exists a 
(Hermitian) deformation quantization of $E$ corresponding to $\qA$,
which is unique up to equivalence. Moreover, if $\star$ is
local/differential/of Vey type, the deformation of $E$ can be chosen
to be of the same type. 
\end{theorem}
\begin{proof}
Existence and uniqueness of (Hermitian) deformations follow from
Propositions \ref{DefModProp}, \ref{HerEProp} and the observation
before the theorem.

Suppose now $\star$ is local/differential/of Vey type. Choose a 
deformation $\qP$ of $P_0$ and let us consider the 
$\qA$-action on $\mathcal{E}[[\lambda]]$ induced by $I$ as in
(\ref{IEq}). Note that if we write 
$I=\sum_{r=0}^\infty I_r \lambda^r$, it follows from 
(\ref{IrEq}) that each 
$I_r: \mathcal{E} \longrightarrow M_n(\mathcal A)$
is local/differential/differential of order r. 
Moreover, $I^{-1}= \sum_{r=0}^\infty J_r \lambda^r$ has the same
property. From $x \bullet A = I^{-1} ( \qP \star x \star A)$ and
$\qh(x,y) = \SP{\qP\star x, \qP\star y}$, it follows directly that
$\Right_r$ and $h_r$ have the same desired properties.
\end{proof}

It is well-known that the complex algebras 
$\mathcal A = C^\infty(M)$ and 
$\mathcal B = \Gamma^\infty(\End(E))$ are Morita equivalent and hence,
as such, have the same algebraic deformation theory (see
\cite[Sect.~16]{GS88}). We will now observe that these algebras in fact
have the same \emph{local} and \emph{Hermitian} deformation theories.
See \cite{LecRog88} for the case of local deformations of a trivial bundle.

Let $\Mloc(\mathcal B)=\{\mbox{local deformations of } \mathcal B\}$,
i.e. elements in $\Mloc(\mathcal B)$ are local star-products $\star$
on $\mathcal{B}[[\lambda]]$. Suppose 
$\star_1, \star_2 \in \Mloc(\mathcal B)$ are equivalent, through an
equivalence transformation
$\qT = \id + \sum_{r=1}^\infty T_r \lambda^r$. We call
$\qT$ \emph{local} if each $T_r$ is local. 
We remark that the argument in \cite[Lem.~1.1.4]{Deligne95} shows that any
equivalence transformation between local star-products is
automatically local. We define $\Def(\mathcal B)$ as the quotient of
$\Mloc(\mathcal B)$ by (local) equivalences. Similarly, we define
$\Def^*(\mathcal B)$ to be the set of local Hermitian deformations of
$\mathcal B$, up to equivalence. We denote the equivalence class of
$\star$ in $\Def^{(*)}(\mathcal B)$ by $[\star]$.


Recall that if we identify $\Gamma^\infty(E)=P_0\mathcal{A}^n$, 
where $P_0 \in M_n(\mathcal A)$ is an idempotent (projection, in the
Hermitian case), we can
write $\mathcal B = P_0M_n(\mathcal A)P_0$.
Let $\qP \in M_n(\qA)$ be an idempotent (projection)
deforming $P_0$. Note that, by Corollary \ref{PAPCor}, the
algebra $\qP \star M_n(\qA) \star \qP$ defines a (Hermitian) deformation
of $\mathcal B$ via $I$. Moreover, it follows from (\ref{IrEq}) that,
if $\star$ is local (or differential/of Vey type), then so is the
deformation of $\mathcal B$ defined by $I$.
(It would be interesting to compare this construction
of star-products on $\mathcal B$ with Fedosov's construction 
\cite[Sect.~7]{Fed94a} in the symplectic case.) Note also that if $\qP'$
is another deformation of $P_0$, then $\qP' \star M_n(\qA) \star \qP'$
and $\qP \star M_n(\qA) \star \qP$ induce equivalent deformations of
$\mathcal B$. To see that, recall that 
$\qP \star M_n(\qA) \star \qP = \End_{\qA}(\qP \star \qA ^n)$ 
and $\qP' \star M_n(\qA) \star \qP' = \End_{\qA}(\qP' \star \qA ^n)$ 
and $\qP \star \qA ^n$ and $\qP' \star \qA ^n$ are equivalent as
$\qA$-modules by Proposition~\ref{DefModProp}. It is simple to check
that this gives rise to a well-defined map 
$\Phi: \Def^{(*)}(\mathcal A) \longrightarrow \Def^{(*)}(\mathcal B)$.
\begin{proposition}
\label{BijecProp}
The map $\Phi: \Def^{(*)}(C^\infty(M)) 
\longrightarrow \Def^{(*)}(\Gamma^\infty(\End(E)))$ is a bijection.
\end{proposition}
\begin{proof}
As we have remarked, the algebras $\mathcal A = C^\infty(M)$ and
$\mathcal B = \Gamma^\infty(\End(E))$ are Morita equivalent,
and $\mathcal E = \Gamma^\infty(E) = P_0 \mathcal{A}^n$ is a $(\mathcal B$-$
\mathcal A)$-bimodule defining this equivalence. 
By symmetry of Morita equivalence, it follows that there exists an idempotent
$Q_0 \in M_m(\mathcal B)$, for some $m \geq 1$,
so that $\mathcal E = \mathcal{B}^mQ_0$ ($\mathcal B^m$ as \emph{row} vectors)
as a left $\mathcal B$-module and 
$\mathcal A = \End_{\mathcal B}(\mathcal{B}^mQ_0) = Q_0M_n(\mathcal B)Q_0$ 
as a unital algebra. In the 
Hermitian case, we note that, by \cite[Thm.~26]{Kapla68}, we can 
actually choose $Q_0$ satisfying $Q_0^* = Q_0$. Note that since $\mathcal A$
is commutative (and so is $Q_0M_n(\mathcal B)Q_0$), it follows from
\cite[Cor.~7.7]{BuWa99a}, \cite[Thm.~4.2]{Ara99a} that in fact $\mathcal A$ and
$Q_0M_n(\mathcal B)Q_0$ are isomorphic as $^*$-algebras (these algebras
are Morita $^*$-equivalent in the sense of \cite{Ara99a}). Therefore, we can define
a map $\hat{\Phi}: \Def^{(*)}(\mathcal B) \longrightarrow \Def^{(*)}(\mathcal A)$
just as we did for $\Phi$.

Let $[\star] \in \Def^{(*)}(\mathcal A)$ and $\qA = (\mathcal{A}[[\lambda]],
\star)$. Let $\qP$ be a deformation of $P_0$ and $\qB=(\mathcal{B}[[\lambda]],
\star')$ be the deformation induced by $\qP\star M_n(\qA)\star \qP$, in such
a way that $[\star']=\Phi([\star])$. Let $\qE=\qP\star\qA^n$, which is
naturally a $(\qB$-$\qA)$-bimodule.
Note that, by Morita theory, we have $\qA = \End_{\qB}(\qE)$.
Now pick $\qQ \in M_n(\qB)$, a deformation
of $Q_0$ and let $\qE'=\qB^m \star' \qQ$.
Then $[\star'']=\hat{\Phi}([\star'])=\hat{\Phi}\circ\Phi([\star])$
is induced by $\qQ \star' M_m(\qB) \star' \qQ$, and 
$\qA' = (\mathcal{A}[[\lambda]], \star'' )$ can be identified with
$\End_{\qB}(\qE')$. 
 
Finally, it is not hard to check (see the existence part of Proposition
\ref{DefModProp}) that $\qE$ and $\qE'$ are both  deformations of 
$\mathcal E = P_0\mathcal{A}^n =\mathcal{B}^m Q_0$ corresponding to
$\qB$. It then follows from Proposition \ref{DefModProp} (for left
modules) that these deformations are equivalent and hence so are
$\star$ and $\star''$. Therefore $\hat{\Phi} \circ \Phi = \id$ and a
similar argument shows that $\Phi \circ \hat{\Phi} = \id$. This
concludes the proof. 
\end{proof}

\subsection{The semi-classical limit}
\label{SemiClassSec}

We shall now compute the first order term of the deformed
module $\qE = (\mathcal E[[\lambda]], \bullet)$ over 
$\qA = (\mathcal A[[\lambda]], \star)$, where $\mathcal A$ is a unital
$^*$-algebra, $\star$ a Hermitian deformation and 
$\mathcal E = P_0 \mathcal A^n$ for $P_0 \in M_n(\mathcal A)$ a projection. 
The deformed module structure $\bullet$ is defined via a deformation 
$\qP \in M_n (\qA)$ of $P_0$ and the isomorphism 
$I: P_0 \mathcal A^n [[\lambda]] \to \qP \star \qA^n$ as in
Lemma \ref{ILem}. 
A simple computation yields
\begin{equation}
    \label{RoneTerm}
    R_1 (x, A) = P_0 C_1 (x, A)
\end{equation}
for the first order term of $\bullet$, where $C_1 (x, A)$ is defined
by $C_1 (x, A)_i = C_1 (x_i, A)$ for $x \in \mathcal A^n$ and 
$A \in \mathcal A$. In particular $R_1$ does not depend on the chosen
deformation of $P_0$.

Let us assume that $\mathcal A$ is commutative. It is well-known
that in this case the skew-symmetric part of $C_1$ is a Poisson
bracket for $\mathcal A$. In order to get a \emph{real} Poisson
bracket, we use the convention
\begin{equation}
    \label{PoissonBracket}
    \{A_1,A_2\} := 
    \frac{1}{\im}(C_1(A_1,A_2) - C_1(A_2, A_1)), 
\end{equation}
for $A_1, A_2 \in \mathcal A$. Thus 
$\{A_1, A_2\}^* = \{ A_1^*, A_2^*\}$ follows from the fact that
$\star$ is a Hermitian deformation.
Let us assume furthermore that $C_1$ is skew-symmetric (this yelds no
loss of generality for star-products since any differential cocycle $C_1$ 
in this case is cohomologous to its skew-symmetric part \cite{Vey75}).
Then $C_1(A_1, A_2)= \frac{\im}{2}\{A_1,A_2\}$.
Let us consider the bracket  
$\{\cdot,\cdot\}_{\mathcal E}: 
\mathcal{E} \times \mathcal A \longrightarrow \mathcal{E}$ given by
\begin{equation}
    \label{BrkEEq}
    \{x,A\}_{\mathcal{E}} 
    = P_0\{x,A\} 
    = \frac{2}{\im}R_1(x,A), \;
    x\in \mathcal{E}, A \in 
    \mathcal A
\end{equation}
obtained from the semi-classical limit of $\bullet$.
\begin{proposition}
    \label{PoissonModuleProp}
    The bracket $\{\cdot,\cdot\}_{\mathcal E}$ defines the structure
    of a right Poisson module in the sense of \cite[Def.~3.1]{RVW96}
    on $\mathcal E$. It is a Poisson module in the sense of 
    \cite[Def.~3.2]{RVW96} if and only if the curvature 
    $R_{\mathcal E}: 
    \mathcal A \times \mathcal A \to \End_{\mathcal A}(\mathcal E)$,
    defined by
    \begin{equation}
        \label{CurvDef}
        R_{\mathcal E} (A_1, A_2) x =
        \{ \{x, A_1 \}_{\mathcal E}, A_2  \}_{\mathcal E} -
        \{ \{ x, A_2 \}_{\mathcal E},A_1 \}_{\mathcal E} - 
        \{ x, \{A_1,A_2\}\}_{\mathcal E},
    \end{equation}
    vanishes.
\end{proposition}
\begin{proof}
From (\ref{BrkEEq}) it easily follows that
$\{\cdot,\cdot\}_{\mathcal E}$
satisfies the natural Leibniz rules (see the last two equations
in \cite[Eq.~(16)]{RVW96}) and this implies the first statement.
A simple computation shows that
$R_{\mathcal E}(A_1, A_2) \in \End_{\mathcal A} (\mathcal E)$. 
Note that $R_{\mathcal E} = 0$ is equivalent to the first equation in
\cite[Eq.~(16)]{RVW96} and the last statement follows directly from 
\cite[Def.~3.2]{RVW96}
\end{proof}

Suppose $\mathcal A = C^\infty(M)$, where $M$ is a Poisson manifold.
Let $P_0 \in M_n(\mathcal A)$ be a projection
and let $E \to M$ be given by the image of $P_0$, so that 
$\Gamma^\infty(E) = P_0 \mathcal{A}^n$. Let $d$ denote the natural flat
connection defined on the trivial vector bundle 
$M \times \mathbb{C}^n \to M$, given by component-wise exterior
differentiation. Then $\nabla := P_0 \circ d$ defines a connection
on $E$, sometimes called the \emph{Levi-Civita} connection of $E$.

\begin{corollary}
Let $\mathcal E = \Gamma^\infty(E)$ and $\{\cdot,\cdot\}_{\mathcal E}$
be as in (\ref{BrkEEq}). Then 
$(\mathcal E, \{\cdot,\cdot\}_{\mathcal E})$ is a Poisson
module (in which case $E$ is called a Poisson vector bundle)
if and only if the Levi-Civita connection $\nabla$ is flat on each
symplectic leaf of $M$.
\end{corollary}
\begin{proof}
Note that if $x \in \Gamma^\infty(E)$ and $f \in C^\infty(M)$, we have
$\{ x, f\}_{\mathcal E}= \nabla_{X_f}x$,
where $X_f$ is the Hamiltonian vector field corresponding to $f$. 
Observe that the curvature tensor corresponding to $\nabla$, $R^{\nabla}$,
satisfies
$R_{\mathcal E}(f,g)x = R^\nabla(X_f, X_g)x$,
for all $f,g \in C^\infty(M)$ and $x \in \Gamma^\infty(E)$.
This implies the result.
\end{proof}

We note that $\{\cdot,\cdot\}_{\mathcal E}$  defines a 
linear contravariant connection on $E \to M$ (by $D_{df}x=
\{x,f\}_{\qE}$), which is just the one induced by
$\nabla$ \cite[Sect.~2]{Rui2000}.

For a Hermitian star-product $\star$ on $M$, let $\star'$ be the
corresponding deformation of 
$\Gamma^\infty(\End(E)) = P_0 M_n(\mathcal A) P_0$ 
induced by a deformation $\qP$ of $P_0$ and
the isomorphism 
$I: P_0M_n(\mathcal A)P_0[[\lambda]] \longrightarrow \qP 
\star M_n(\qA) \star \qP$
as in Corollary \ref{PAPCor}. 
Note that the center of $\Gamma^\infty(\End(E))$, denoted by $Z$,
is isomorphic to $C^\infty(M)$ through $f \mapsto f P_0$.
As discussed in \cite[Prop.~1.2]{RVW96}, the skew-symmetric part of
the semi-classical limit of
$\star'$ endows the pair $(\Gamma^\infty(\End(E)), Z)$ with the 
structure of a \emph{Poisson fibred algebra} (see \cite[Def.~1.1]{RVW96}).
For $L,S \in M_n(\mathcal A)$, 
let $\{ L, S \} = \frac{1}{\im}(C_1(L,S) - C_1(S,L))$, where
$C_1(L,S) \in M_n(\mathcal A)$ is defined by
${C_1(L,S)}_{i,j}= \sum_{r=1}^nC_1(L_{i,r},S_{r,j})$.

\begin{proposition}
The Poisson fibred algebra bracket $\{\cdot, \cdot \}':
\Gamma^\infty(\End(E)) \times Z
\longrightarrow \Gamma^\infty(\End(E))$ induced by
$(\Gamma^\infty(\End(E))[[\lambda]],\star')$ is given by
$\{L_0, u\}' = P_0 \{ L_0, u \} P_0$, 
for $L_0 \in \Gamma^\infty(\End(E))$ and $u \in Z$. Moreover, the Poisson
bracket defined by $\{\cdot, \cdot \}'$ on $Z$ coincides
with the original Poisson bracket $\{\cdot, \cdot \}$ on $C^\infty(M)$.
\end{proposition}

\begin{proof}
For $L_0, S_0 \in P_0M_n(\mathcal A)P_0 = \Gamma^\infty(\End(E))$,  write
$ L_0 \star' S_0 = I^{-1}(I(L_0)\star I(S_0)) = \sum_{r=0}^\infty B_r(L_0, S_0)
\lambda^r$. It is not hard to check that
\begin{equation}
B_1(L_0, S_0)= P_0 C_1(L_0, S_0)P_0.
\end{equation}
It is then clear that 
$\{L_0,S_0\}':= \frac{1}{\im}(B_1(L_0,S_0) - B_1(S_0, L_0))
 = P_0\{L_0, S_0\} P_0$.
 
The bracket $\{\cdot, \cdot \}'$ defines an action of the Poisson algebra
$(Z,\{\cdot, \cdot \}' )$ on $\Gamma^\infty(\End(E))$ by derivations and this
implies that $P_0\{P_0, \cdot \} P_0 = P_0\{\cdot, P_0 \} P_0 = 0$.
Hence the Leibniz rule for $\{ \cdot, \cdot \}$ yields
$ \{ f P_0, g P_0 \}' = P_0 \{ f, g \} P_0 = \{f,g\}P_0$, 
$f, g \in C^\infty(M)$, and this immediately
shows that the bracket $\{ \cdot, \cdot \}'$ on $Z$ coincides with
$\{ \cdot, \cdot \}$ after the identification $Z \cong C^\infty(M)$. 
\end{proof}

\section{Strongly full projections and formal Morita equivalence} 
\label{FullProjSec}

It is known that unital algebras which are Morita equivalent have equivalent
algebraic deformation theory and, moreover, corresponding deformations are again
Morita equivalent (see \cite[Sect.~16]{GS88}). 
We will show in this section that local, Hermitian 
deformations which are related by $\Phi$ are actually formally Morita
equivalent, which is a notion 
stronger than the classical Morita equivalence and related to strong Morita
equivalence of $C^*$-algebras (see \cite{Rief82}). We now briefly recall
the definitions, see \cite{BuWa99a,BuWa2000} for details.

Let $\mathcal A, \mathcal B$ be
$^*$-algebras over $\ring C = \ring R(\im)$,
where $\ring R$ is an ordered ring. Consider a 
$(\mathcal B$-$\mathcal A)$-bimodule $\mathcal E$ with an 
$\mathcal A$-valued, $\mathcal A$-right linear positive semi-definite,
full inner product $\SP{\cdot,\cdot}$ as well as a 
$\mathcal B$-valued, $\mathcal B$-left linear, positive semi-definite,
full inner product $\Theta_{\cdot,\cdot}$ such that 
$\SP{x, By} = \SP{B^*x, y}$, $\Theta_{x,yA} = \Theta_{xA^*, y}$, and
$\Theta_{x,y} z = x \SP{y, z}$ for all $A \in \mathcal A$, 
$B \in \mathcal B$, and $x, y, z \in \mathcal E$. Here positivity is
understood in the sense of positive algebra elements of a
$^*$-algebra (see Section~\ref{StarAlgSect}). Fullness means 
$\SP{\mathcal E, \mathcal E} = \mathcal A$ and 
$\Theta_{\mathcal E, \mathcal E} = \mathcal B$, respectively. Let
$(\pi,\mathfrak H)$ be a $^*$-representation of $\mathcal A$ on a
pre-Hilbert space over $\ring C$ and consider the $\mathcal
A$-balanced tensor product 
$\mathcal E \otimes_{\mathcal A} \mathfrak H$ with the inner product 
$\SP{x \otimes \psi, y \otimes \phi} = \SP{\psi, \pi(\SP{x,y})\phi}$.
If this inner product is positive semi-definite for all 
$(\pi, \mathfrak H)$ then $\mathcal E$ is said to satisfy property
(P). The analogous property for $^*$-representations of
$\mathcal B$ is called (Q). If $\mathcal E$ satisfies all
these requirements then $\mathcal E$ is called an \emph{equivalence
bimodule} and $\mathcal A$ and $\mathcal B$ are called \emph{formally Morita
equivalent}. If in addition the actions of $\mathcal A$ and 
$\mathcal B$ on $\mathcal E$ are non-degenerate then $\mathcal E$ is
called non-degenerate. Note that we are dealing here with unital
$^*$-algebras only. We remark that this purely algebraic notion is
equivalent to strong Morita equivalence when applied to
$C^*$-algebras, see \cite{BuWa2000}.

Let $\mathcal A$ be a $^*$-algebra over $\ring C$ and 
$P_0 \in M_n(\mathcal A)$ be a projection. Let 
$\mathcal E = P_0 \mathcal{A}^n$, considered as a right
$\mathcal A$-module and left
$\End_{\mathcal A}(\mathcal E)= P_0M_n(\mathcal A)P_0$ module. 
Note that both actions are non-degenerate.
We noted in Section \ref{InnerProdSec} that $\mathcal E$
has a canonical $\mathcal A$-valued inner product $h_0$.
Consider the map $\Theta_{\cdot, \cdot}: \mathcal E \times \mathcal E
\longrightarrow \End_{\mathcal A}(\mathcal E)$, defined by
$\Theta_{x,y}z = x h_0(y,z)$, for $x,y,z \in \mathcal E$.
The following definition will give a sufficient condition to guarantee
that $\mathcal E$, equipped with $h_0$ and $\Theta_{\cdot,\cdot}$, is an
equivalence bimodule.
\begin{definition}
\label{StrongFullDef}
A projection $P_0 \in M_n (\mathcal A)$ is called strongly full if there
exists an invertible element $\tau \in \mathcal A$ such that 
$\tr P_0 = (\tau \tau^*)^{-1}$. 
\end{definition}
It turns out that this is indeed a stronger version of the usual
notion of full projections (recall that $P_0 \in M_n(\mathcal A)$
is full if $M_n(\mathcal A)P_0M_n(\mathcal A)=M_n(\mathcal A)$).
\begin{theorem}
\label{AlmostMainTheoremII}
Let $P_0 \in M_n (\mathcal A)$ be a strongly full projection. Then
$\mathcal E$ is a non-degenerate equivalence bimodule and thus
$\End_{\mathcal A} (\mathcal E)$ and $\mathcal A$ are formally Morita
equivalent.
\end{theorem}
\begin{proof}
Let $x, y \in \mathcal E \subseteq \mathcal A^n$. Then a
straightforward computation shows the relation
\begin{equation}
\label{NiceI}
    \sum_i \Theta_{x, P_0 e_i \tau} \Theta_{P_0 e_i \tau, y} =
    \Theta_{x, y}.
\end{equation}
It follows immediately that $\Theta_{x, x}$ is a positive algebra
element and that $\Theta_{\cdot,\cdot}$ is full. The fullness of
$\SP{\cdot,\cdot}$ follows from
\begin{equation}
\label{NiceII}
    \sum_i \SP{P_0e_i\tau, P_0e_i\tau} = \tau^* \tr P_0 \tau = \Unit.
\end{equation}

We observe that property (P) can be easily shown as in 
\cite[Sect.~6]{BuWa99a}.
In order to prove (Q), we recall that this just means
(P) for the complex-conjugate bimodule $\cc{\mathcal E}$.
Equivalently, we can consider $^*$-representations of 
$\End_{\mathcal A} (\mathcal E)$ on pre-Hilbert
spaces from the \emph{right}. Thus let $(\pi, \mathfrak H)$ be such a
$^*$-representation of $\End_{\mathcal A}(\mathcal E)$ from the right
and let $\phi_1, \ldots, \phi_r \in \mathfrak H$ as well as 
$x_1, \ldots, x_r \in \mathcal E$. Then by (\ref{NiceI})
\[
    \sum_{i,j} \SP{\phi_i \otimes x_i, \phi_j \otimes x_j}
    =
    \sum_{i,j} \SP{\phi_i \pi\left(\Theta_{x_i, x_j}\right), \phi_j}
    =
    \sum_{i,j,k} \SP{\phi_i \pi\left(\Theta_{x_i, P_0e_k\tau}\right),
                     \phi_j \pi\left(\Theta_{x_j, P_0e_k\tau}\right)}
    \ge 0
\]
and this
shows the positivity needed for (Q). This concludes the proof.
\end{proof}

We remark that, by Lemma~\ref{TrickLem}, it immediately follows that 
deformations of strongly full projections are again strongly full. 
\begin{proposition}
\label{StrongFullDeformProp}
Let $\qA$ be a Hermitian deformation of a $^*$-algebra $\mathcal A$
and $P_0 \in M_n (\mathcal A)$ be a strongly full projection. 
Then every deformation $\qP \in M_n (\qA)$ of $P_0$ is again strongly full 
and therefore
$\End_{\mathcal \qA} (\qP \star \qA^n)$ and $\qA$ are formally Morita
equivalent via the non-degenerate equivalence bimodule $\qE = \qP \star \qA^n$.
\end{proposition}

Note that if $\mathcal A = C^\infty(M)$ and $P_0 \in M_n(\mathcal A)$
is a nowhere zero projection, then $\tr P_0$ is nonzero and constant
on connected components of $M$.
Hence $P_0$ is strongly full. 
\begin{theorem}
\label{MainTheoremII}
Let $E \to M$ be a (non-zero) Hermitian vector bundle over a Poisson
manifold $M$ with star-product $\star$. Then every deformation $\qE$
of $\mathcal E = \Gamma^\infty (E)$ is a non-degenerate equivalence
bimodule and therefore $\End_{\qA} (\qE)$ and $\qA$ are formally
Morita equivalent.
\end{theorem}
\begin{corollary}
Consider the bijective map $\Phi:\Def^*(C^\infty(M)) \longrightarrow
\Def^*(\Gamma^\infty(\End(E)))$ as in Proposition
\ref{BijecProp}. If $\qA = (C^\infty(M)[[\lambda]], \star)$ and
$\qB = (\Gamma^\infty(\End(E))[[\lambda]], \star')$ are Hermitian
deformations such that $[\star']=\Phi([\star])$, then $\qA$ and $\qB$
are formally Morita equivalent. 
\end{corollary}

The formal Morita equivalence of $C^\infty(M)$ and 
$\Gamma^\infty(\End(E))$ follows from Theorem \ref{MainTheoremII} by
considering the
undeformed product for the trivial Poisson bracket. This was shown
in \cite[Sect.~6]{BuWa2000} by a more direct argument.
We also remark that formal Morita equivalence implies that the
$^*$-representation theories of the involved $^*$-algebras on
pre-Hilbert spaces over $\mathbb{C}[[\lambda]]$ are the
same,
see~\cite[Thm.~5.10]{BuWa99a}.

%
%

\section*{Acknowledgements}

We would like to thank Michael Artin, Martin Bordemann, 
Michel Cahen, Laura DeMarco,
Rui L. Fernandes, Viktor Ginzburg,
Simone Gutt, Tim Swift and Alan Weinstein for valuable
discussions and remarks. 

%
%

\begin{small}

\end{small}
\end{document}